\documentclass[twoside,11pt]{article}

\usepackage{amsfonts}
\usepackage{amssymb}
\usepackage{amsmath}
\usepackage{exscale,cmmib57}
\usepackage{authblk}
\usepackage{fancyhdr}
\usepackage{mdframed}
\usepackage[colorlinks = true,
linkcolor = red,
urlcolor  = blue,
citecolor = blue,
anchorcolor = blue]{hyperref}
\usepackage{stmaryrd}

\newtheorem{thm}{\bfseries Theorem}

\newtheorem{lem}{\bfseries Lemma}

\newtheorem{cor}{\bfseries Corollary}
\newtheorem{rem}{\it Remark}

\newcommand{\R}{\mathbb{R}}

\newcommand{\dom}{\mathop{\rm dom}}
\newcommand{\Int}{\mathop{\rm int}}

\newcommand{\Min}{\mathop{\rm Min}}

\newcommand{\argmin}{\mathop{\rm argmin}}

\title{\large\bf In the blessed memory of Prof. Dr. Hedy Attouch\\[.5cm]
\Large\bf $\varepsilon$-OPTIMALITY IN REVERSE CONVEX OPTIMIZATION}
\author{\sc M. El Maghri\footnote{Corresponding author. E-mail: elmaghri@yahoo.com. ORCID: 0000-0001-6642-0145.} and H. Sellak}
\affil{\small\em Department of Mathematics and Computer, Faculty of Sciences A\"{\i}n Chock, Hassan II University, 
	Casablanca, BP. 5366, Morocco}
\date{\today}

\begin{document}

\maketitle

\vspace{-0.2cm}
\noindent
{\small

\noindent
{\sl Abstract.} 
We characterize approximate global optimal solutions ($\varepsilon$-optima) to reverse 
optimization problems, namely, problems whose non-convex constraint is of the form 
$h(x)\geq 0$. This issue has not been addressed previously in the literature. Our idea 
consists of converting the reverse program into an unconstrained bicriteria DC program. 
The main condition presented is obtained in terms of Fenchel's $\varepsilon$-subdifferentials 
thanks to an earlier result in difference vector optimization by El Maghri. This extends 
and improves similar results from the literature dealing with exact ($\varepsilon=0$) solutions. 
Moreover, as we consider functions with extended values, our approach also applies to reverse 
problems subject to additional convex constraints, provided that Moreau--Rockafellar or 
Attouch--Brézis constraint qualification conditions are satisfied. Similarly, new results  
for the special case of a nonlinear equality constraint $h(x) = 0$ are also obtained. \\[.3cm] 
}

\vspace{0.2cm}

\section{Introduction}

Nonconvex minimization has long been a difficult global optimization problem 
that has attracted a growing interest in recent years. Motivation mainly comes 
from modelling real-world applications in applied mathematics and operational 
research (see e.g., \cite{hor}), where one may be confronted with a 
reverse optimization problem of the form\\
$$
\begin{array}{ll}
\mbox{(ROP)}\quad&\Min\ f(x)\\[.1cm]
&\left\{\begin{array}{l}
h(x)\geq 0,\\[.1cm]
x\in C
\end{array}\right.
\end{array}
$$ 
where $f$ and $h$ are convex functions over a convex subset $C$ of a real topological 
vector space $X$, while $h$ is assumed to be a nonconcave function. 

\smallskip 

The main difficulty for dealing with (ROP) is due to the reverse inequality constraint 
$h(x)\geq 0$, which can generate nonconvex or even disconnected feasible regions, so that 
local minima may not be global. Such constraints first appeared in Rosen's paper \cite{ros} 
as part of a discretization process in optimal control, and subsequently, in several papers 
covering many application areas, such as engineering design \cite{avr}, networks \cite{ban}, 
economics \cite{hil}, production design \cite{vid}, management \cite{zale}, and so on (see 
e.g., references given therein). 

\smallskip

Reverse optimization has been widely studied from a numerical point of view. The main 
approaches proposed, basically, either use techniques of local optimization focusing on 
stationary points generally of KKT type (see e.g., \cite{avr,ros,vid}) or proceed with 
global optimization strategies, such as branch-and-bounds, branch-and-cuts, outer 
approximation, concave programming, etc (see e.g., \cite{ban,ben,hil,muu,tuy,yam,zha}).

\smallskip 

On the other hand, important contributions have been devoted to the theoretical aspects 
of global optimal solutions and some interesting properties of reverse problems like (ROP). 
For instance, Tuy \cite[Proposition 5.1]{tuy} established, under a stability condition, 
a duality result between (ROP) and a convex maximization problem (or equivalently, a 
concave minimization problem) assuming that $f$ and $h$ are finite convex functions. Later, 
Strekalovsky \cite[Theorem 4]{str} developed a characteristic condition of global optimality 
given in terms of subdifferentials. An analogous characterization was proposed by Hiriart-Urruty 
\cite[Theorem 3.5]{hir} in terms of $\varepsilon$-subdifferentials, where $f$ and $h$ are assumed 
to be finite convex functions. Likewise, the paper by Tseveendorj \cite{tse} provided a weaker 
characteristic condition of the type given in \cite{str}, but also restricts $f$ and $h$ to 
finite convex functions. However, the issue of characterizing {\em approximate} optimal solutions 
(or $\varepsilon$-optima) for (ROP) has not yet been addressed in the literature (even when 
$C=X=\mathbb{R}^n$). In this paper, we attempt to provide answers to this question.
 
\smallskip 

From our perspective, it will be interesting to recall that many global optimization problems 
can be converted into reverse optimization problems as previously discussed in several papers. 
For instance, the binary constraint $x_i=0$ or $1$ that arises in integer programming can be 
reviewed as a reverse constraint $x_i^2 - x_i\geq 0$ and  the box-constraint $0\leq x_i\leq 1$.  
Also, the classic condition of complementarity of the form  $\sum_i x_iy_i = 0$, with $x_i, 
y_i\geq 0$, can be transformed into the reverse convex constraint $\sum_i\max(-x_i,-y_i)\geq 0$. 
In turn, the problem of minimizing a concave function $g$ itself is obviously convertible to the 
problem of minimizing an artificial variable $t$ under the reverse convex constraint $t-g(x)\geq 0$. 

\smallskip 

In this paper, we show that the reverse constraint $h(x)\geq 0$ can always be penalized in 
such a way that (ROP) becomes equivalent, in a certain sense, to an unconstrained bicriteria 
DC optimization problem. This allows us, via an earlier result by El Maghri \cite{elm1} in 
difference vector optimization, to derive an approximate optimality condition for (ROP) that 
recovers the result given in \cite{hir} for exact solutions. On the other hand, since the 
functions involved in our context are assumed to have extended values, our approach also applies 
to the case where the reverse constraint coexists with other convex constraints, provided that 
Moreau--Rockafellar's qualification condition or the weaker condition of Attouch--Brézis is 
satisfied. Another interesting application concerns the special case of the nonlinear equality 
constraint $h(x)= 0$. As we show, this constraint can be viewed as a reverse constraint with 
an additional convex constraint, or it can be seen as a single reverse convex constraint provided 
$f$ and $h$ are finite convex functions. For both cases, new characteristic conditions of approximate 
optimal solutions are given.   

\medskip 

\section{A bicriteria approach}

Let us first consider the following general vector optimization problem:\\[.1cm]
$$
\mbox{(VOP)}\qquad \epsilon\text{\bf-}\Min_{x\in S}\ F(x)
$$\\
where $F : X \supseteq S\to Y\sqcup\{+\infty\}$, $X$ and $Y$ are real topological 
vector spaces with $Y$ separated. It is adjoined to $Y$ the abstract maximal element 
$+\infty$ obeying the operations $y \pm (+\infty)=+\infty\;$ ($\forall\: y\in Y\sqcup\{+\infty\}$), 
$\alpha\cdot (+\infty)=+\infty$ if $\alpha\in\R_+\setminus\{0\}$ and  $0\cdot(+\infty)=0$. 
Also the space $Y$ is assumed to be endowed with a nonempty convex cone $Y_+$, with nonempty 
interior ($\Int Y_+\neq\emptyset)$, inducing the preorder relations:\\[.1cm]
$$
\begin{array}{lcl}
y\leq_{Y_+} y' &\Leftrightarrow& y'-y  \in Y_+,\\[.2cm]
y<_{Y_+} y' &\Leftrightarrow& y'-y  \in \Int Y_+,\\[.2cm]
y \lneq_{Y_+} y' &\Leftrightarrow& y'-y\in Y_+\setminus l(Y_+),
\end{array}
$$\\[.1cm]
where $l(Y_+)=Y_+\sqcap -Y_+$ is the lineality of $Y_+$. Recall that when $l(Y_+)=\{0\}$, 
the cone $Y_+$ is said to be pointed and the preorder induced by $Y_+$ becomes an order. 
 
\medskip

Given $\epsilon\in Y$, the sets of $\epsilon$-$\sigma$-efficient points for (VOP) with 
respect to the preorder cone $Y_+$, are the sets of the optimal $\epsilon$-solutions 
to (VOP) taken either in the strong Pareto sense, or in the efficient Pareto sense, or in 
the weak Pareto sense, or in the proper Pareto sense depending on the choice of 
$\sigma\in\{s,e,w,p\}$. These sets are given under the following unified form: 
for $\sigma\in\{s,e,w\}$, \\[.1cm]
$$
\begin{array}{c}
E_\epsilon^\sigma(F,\;S,\;Y_+)=\{\bar{x} \in S\sqcap\dom F :\; \forall x\in S,\;  
F(x)\not<_{Y_+}^\sigma F(\bar{x}) - \epsilon\},\\[.5cm]
E_\epsilon^p(F,\;S,\;Y_+)= \displaystyle{\bigsqcup_{\widehat{Y}_+\in 
		\mathcal{C}(Y_+)}E_\epsilon^e(F,\;S,\;\widehat{Y}_+)},
\end{array}
$$\\
where $\mathcal{C}(Y_+):=\{\widehat{Y}_+\subsetneq Y\; \mbox{convex cone} :  Y_+\setminus l(Y_+) \subseteq \Int \widehat{Y}_+\}$, 
the effective domain of $F$ is denoted by $\dom F = \{x\in X \ :\; F(x)\in Y\}$, and, the symbol 
``$\not<_{Y_+}^\sigma$'' stands for the opposite relation to:\\
$$
y<_{Y_+}^\sigma y' \quad\Leftrightarrow\quad \left\{\begin{array}{lcl}
y\not\geq_{Y_+} y'&\quad&\mbox{if } \sigma=s,\\[.2cm]
y<_{Y_+} y'&&\mbox{if } \sigma=w,\\[.2cm]
y\lneq_{Y_+} y'&\quad&\mbox{if } \sigma\in\{e,p\}.
\end{array}\right.
$$

\vskip.3cm
\noindent 
The other opposite relations are denoted like $\not\leq_{Y_+}$, $\not<_{Y_+}$ and $\not\lneq_{Y_+}$, 
reverse ones like $>_{Y_+}^\sigma$, $\geq_{Y_+}$, $>_{Y_+}$ and $\gneq_{Y_+}$.

\bigskip

It is worth noting that (see \cite{elm3} for more details) \\
$$
\begin{array}{c}
E_\epsilon^\sigma(F,S,Y_+)\neq\emptyset\quad \Longrightarrow\quad \epsilon\not<_{Y_+}^\sigma 0.
\end{array}
$$ 

\vskip.3cm

On the other hand, it is immediate that \\
$$
E_\epsilon^s(F,\;S,\;Y_+)\subseteq E_\epsilon^p(F,\;S,\;Y_+)\subseteq 
E_\epsilon^e(F,\;S,\;Y_+)\subseteq E_\epsilon^w(F,\;S,\;Y_+).
$$

\medskip 

For instance, the set of usual optimal $\epsilon$-solutions for (VOP) called strongly 
$\epsilon$-efficient set in Pareto language is $E_\epsilon^s(F,S,Y_+) = 
\{\bar{x}\in S\sqcap\dom F \ : \;\forall x\in S,\; F(x)\geq_{Y_+} F(\bar{x})-\epsilon\}$. 
The $0$-$\sigma$-efficiency reduces to (exact) $\sigma$-efficiency: 
$E^\sigma(F,S,Y_+)=E_0^\sigma(F,S,Y_+)$. In the scalar case ($Y=\R$), $\epsilon$-$\sigma$-efficiency 
reduces to $\epsilon$-suboptimality: \\
\begin{equation}\label{eopta}
E_\epsilon^\sigma(F,S,\R_+) = \epsilon\text{\bf-}\argmin_{x\in S}F(x) = 
\{\bar{x}\in S\sqcap\dom F\ :\; \inf_{x\in S} F(x)\geq F(\bar{x})-\epsilon\},
\end{equation}

\bigskip

The basic result of this paper relates (ROP) to a bicriteria (VOP) as follows.

\smallskip 

\begin{lem}\label{rb}
Let $f,\, h : X \to \R\sqcup\{+\infty\}$ and $\varepsilon \in \R_+$. Then,  
the following statements holds:\\[.1cm] 
$$
\bar{x} \in \mathbf{\varepsilon}\text{\bf-}\argmin_{h(x)\geq 0} f(x)
\quad \Rightarrow \quad \bar{x} \in E^w_{(\varepsilon,0)}\big((f,-h),X,\R^{2}_{+}\big).
$$ \\[-.2cm]
Conversely, if $h(\bar{x})=0$,\\[-.1cm]
$$
\bar{x} \in E^e_{(\varepsilon,0)}\big((f,-h),X,\R^{2}_{+}\big) 
\quad \Rightarrow \quad \bar{x} \in \mathbf{\varepsilon}\text{\bf-}\argmin_{h(x)\geq 0} f(x).
$$
\end{lem}

\noindent
{\em Proof.}\quad Let $\bar{x} \in \mathbf{\varepsilon}\text{\bf-}\argmin_{h(x)\geq 0} f(x)$ 
and suppose that $\bar{x} \not\in E^w_{(\varepsilon,0)}\big((f,-h),X,\R^{2}_{+}\big)$. 
Then, by the very definition, since  $\bar{x} \in \dom f\sqcap\dom(-h)$,  
$$
\exists x \in X \quad \mbox{s.t.} \quad f(x) < f(\bar{x})-\varepsilon  \;\; \mbox{and} \;\; 
-h(x) < -h(\bar{x}) \leq 0.
$$ \\[-.65cm]
Hence, 
$$
\exists x \in X \quad \mbox{s.t.} \quad h(x)>0  \;\; \mbox{and} \;\;  
f(x) < f(\bar{x})-\varepsilon,
$$ \\[-.3cm]
which contradicts the hypothesis. Conversely, assume that $\bar{x}$ is such that 
$h(\bar{x})=0$ and $\bar{x} \in E^e_{(\varepsilon,0)}\big((f,-h),X,\R^{2}_{+}\big)$, 
but $\bar{x} \not\in\mathbf{\varepsilon}\text{\bf-}\argmin_{h(x)\geq 0} f(x)$. Then, 
$$
\exists x \in X \quad \mbox{s.t.} \quad h(x)\geq 0 = h(\bar{x}) \;\; \mbox{and} \;\; 
f(x) < f(\bar{x})-\varepsilon.
$$ \\[-.65cm]
Hence, 
$$
\exists x \in X \;\; \mbox{s.t.} \;\;
(f(x),-h(x))\lneq_{\mathbb{R}^{2}_{+}} (f(\bar{x}),-h(\bar{x}))-(\varepsilon,0),
$$ \\[-.2cm]  
which contradicts the hypothesis that $\bar{x} \in E^e_{(\varepsilon,0)}\big((f,-h),X,\R^{2}_{+}\big)$.\hfill$\Box$

\bigskip 

Hence, following Lemma \ref{rb}, (ROP) can be converted to the bicriteria difference problem:\\
$$
\mbox{(BOP)}\qquad \Min_{x\in X}\ \big(f(x),0\big)-\big(0,h(x)\big) 
$$

\medskip 

Thanks to an earlier result by El Maghri \cite{elm1} on weak $\epsilon$-efficiency in difference 
vector optimization, we will prove in the sequel that the $\varepsilon$-optimal solutions to (ROP) 
can be completely characterized in terms of $\varepsilon$-subdifferentials by means of its equivalent 
problem (BOP). We also need a similar result for $\epsilon$-efficiency in difference vector optimization 
that we prove in the next section. 

\medskip 

\section{An auxiliary result}

Some elements of vector convex analysis and vector global optimization 
(\cite{elm1,elm2,elm3}) will be recalled. The $\epsilon$-$\sigma$-subdifferential 
of a vector mapping $F : X \to Y\sqcup\{+\infty\}$ at $\bar{x}\in \dom F$, is 
defined in $\sigma$-efficient senses of Pareto, with respect to the choice of 
$\sigma\in\{s,p,e,w\}$, by \\[.1cm]
$$
\partial_\epsilon^\sigma F(\bar{x}) = \{ A\in L(X,Y)\ : \quad 
\bar{x}\in E_\epsilon^\sigma(F-A,X,Y_+)\}
$$\\[-.1cm]
where $L(X,Y)$ is the space of linear continuous operators from $X$ to $Y$. This 
definition comes 
from the following immediate property: \\[.1cm]
$$
\bar{x}\in E_\epsilon^\sigma(F,X,Y_+)\quad\Longleftrightarrow\quad 
0 \in \partial_\epsilon^\sigma F(\bar{x}).
$$

\smallskip

If $\epsilon = 0$, then $\partial_0^\sigma F(\bar{x})$ reduces to the exact 
$\sigma$-subdifferential $\partial^\sigma F(\bar{x})$ (see \cite{elm4} for 
more details about this $\sigma$-subdifferential). As usual, 
$\partial_\epsilon^\sigma F(\bar{x})=\emptyset$ if $\bar{x}\not\in\dom F$. 
Notice that 
$\partial_\epsilon^s F(\bar{x}) = \{A\in L(X,Y)\ :\; \forall x\in X,\; F(x)-F(\bar{x})\geq_{Y_+}A(x-\bar{x})-\epsilon \}$ 
is the so-called strong $\epsilon$-subdifferential which is not other than the ordinary 
(Fenchel) $\epsilon$-subdifferential of convex analysis extended to the vector case. 
The efficient $\epsilon$-subdifferential is in turn formulated by
$\partial_\epsilon^e F(\bar{x}) = 
\{A\in L(X,Y)\ :\; \forall\, x\in X,\; F(x)-F(\bar{x})\not\lneq_{Y_+}A(x-\bar{x})-\epsilon \}$. 
Similar inclusions with the $\epsilon$-$\sigma$-efficient sets follow easily:\\
\begin{equation}\label{isdp}
\partial_\epsilon^s F(\bar{x})\subseteq\partial_\epsilon^p F(\bar{x})\subseteq 
\partial_\epsilon^e F(\bar{x})\subseteq\partial_\epsilon^w F(\bar{x}).
\end{equation}\\[-.7cm]

In the finite-dimensional 
space $Y=\mathbb{R}^r$ with $Y_+=\mathbb{R}^r_+$, it is also easily seen that the 
strong $\epsilon$-subdifferential of $F=(f_1,\ldots,f_r)$, for 
$\epsilon= (\varepsilon_1,\ldots,\varepsilon_r)\geq_{Y_+} 0$,  reduces to \\
\begin{equation} \label{ss}
\partial^s_{\epsilon}F(\bar{x}) = \prod_{i=1}^r \partial_{\varepsilon_i} f_i(\bar{x}).
\end{equation} 
In scalar case ($r=1$), all these sets coincide with the classical 
$\epsilon$-subdifferential ($\partial_\epsilon^s F$) usually denoted by 
$\partial_\epsilon F$. Recall in this case the property:   
$\partial_\epsilon(\lambda F)(x)=\lambda\partial_{\frac{\epsilon}{\lambda}} F(x)$
for all $\lambda \in ]0,+\infty[$.
		
\medskip

The polar cone $Y_+^*$ of $Y_+$ is the set of $\lambda\in Y^*$ such that 
$\lambda(Y_+)\subseteq \R_+$, while the strict polar cone $(Y_+^*)^\circ$ is the 
set of $\lambda\in Y^*$ such that $\lambda(Y_+ \setminus l(Y_+))\subseteq \R_+\setminus\{0\}$. 
Obviously, $(Y_+^*)^\circ\subseteq Y_+^*\setminus\{0\}$. We unify the notation of 
the polar cones of $Y_+$ too by putting \\
$$
Y_+^\sigma = \left\{\begin{array}{lcl}
Y_+^*\setminus\{0\}&\qquad&\mbox{if } \sigma=w,\\[.3cm]
(Y_+^*)^\circ&&\mbox{if } \sigma=p.
\end{array}\right.
$$

\medskip 

For each $\lambda\in Y_+^*\setminus\{0\}$, the scalar function 
$\lambda\circ F : X \to \R\sqcup\{+\infty\}$ is defined by 
$\lambda\circ F(x)=\langle \lambda , F(x)\rangle$ if $x\in \dom F$, $+\infty$ else.  
In particular, $\dom (\lambda\circ F) = \dom F$. \\

The following scalarization results will be used.

\smallskip 

\begin{thm} {\em (\cite{elm3})} \label{ssepw}
Let $F : X \to Y\sqcup\{+\infty\}$. Then, for $\sigma\in\{p,w\}$, 
$\forall\epsilon\not<_{Y_+}^\sigma 0$, $\forall \bar{x} \in X$,\\[.2cm]
$$
\partial_\epsilon^\sigma F(\bar{x})\supseteq\hskip.25cm 
\bigsqcup_{\lambda\in Y_+^\sigma}\hskip.25cm\{A\in L(X,Y) :
\lambda\circ A \in \partial_{\langle\lambda,\epsilon\rangle}(\lambda\circ F)(\bar{x})\},
$$\\[-.1cm]
with equality if $F$ is $Y_+$-convex \big[i.e., 
$F(\alpha x+(1 - \alpha) x')\leq_{Y_+} \alpha F(x)+(1-\alpha) F(x')$ for all 
$\alpha \in[0,1]$ and all $x,\,x'\in X$\big], with $Y_+$ pointed as $\sigma=p$.
\end{thm}

\medskip  

The next auxiliary result gives the characteristic condition for ($\epsilon$-)efficiency 
in difference vector optimization. Its proof follows the same lines as its analogue 
initially established in \cite[Theorem 2]{elm1} for the weak ($\epsilon$-)efficiency 
concept which will also be required in the sequel.

\smallskip 

\begin{thm}\label{cnswe} Let $F,\, G : X \to Y\sqcup\{+\infty\}$ be such that $G$ 
satisfies the following hypothesis for $\sigma\in\{e,w\}$: \\[-.3cm]
$$
({\cal H})\qquad\forall\epsilon >_{Y_+}^\sigma 0,\quad \forall x\in\dom G,\quad 
\partial_\epsilon^s G(x) \neq \emptyset.
$$ \\[-.6cm]
Then, $\forall\epsilon\not<_{Y_+}^\sigma 0$, \\
$$
\bar{x}\in E_\epsilon^\sigma (F-G,X,Y_+) \quad \Leftrightarrow \quad
\forall \epsilon'\geq_{Y_+}0, \quad \partial_{\epsilon'}^s\, G(\bar{x})\subseteq 
\partial_{\epsilon'+\epsilon}^\sigma\,F(\bar{x}).
$$
\end{thm}

\smallskip 
\noindent 
{\em Proof.}\quad To prove the theorem for the case $\sigma=e$, we need the following 
lemma. 

\smallskip

\begin{lem}\label{reee}
For all $\epsilon \in Y$, it holds that \\
$$
E^e_{\epsilon}(F,S,Y_+)= \bigsqcap_{\epsilon' \gneq_{Y_+} 0} E^e_{\epsilon'+\epsilon}(F,S,Y_+).
$$
\end{lem}

\smallskip 
\noindent 
{\em Proof of Lemma \ref{reee}.}\quad We must prove that 
$\bar{x} \not\in E^e_{\epsilon}(F,S,Y_{+})$ iff, $\exists \epsilon' \gneq_{Y_+} 0$,  
$\bar{x} \not\in  E^e_{\epsilon'+\epsilon}(F,S,Y_+)$, i.e., \\[-.3cm]
$$
\exists x\in S,\  F(\bar{x})-F(x)-\epsilon \gneq_{Y_{+}} 0 \quad \Leftrightarrow \quad 
\exists \epsilon' \gneq_{Y_{+}} 0,\ \exists x\in S,\  F(\bar{x})-F(x)-\epsilon\gneq_{Y_{+}}\epsilon', 
$$ \\[-.1cm]
which is true taking for the direct implication 
$\epsilon'=\frac{F(\bar{x})-F(x)-\epsilon}{2}$.\hfill$\Box$

\bigskip 

Now, let us go back to the proof of Theorem \ref{cnswe} by considering the two 
implications separately. 

\bigskip
\noindent
Necessity: suppose that $\partial_{\epsilon'}^s\, G(\bar{x})\not\subseteq\partial_{\epsilon'+\epsilon}^e\,F(\bar{x})$ 
for some $\epsilon'\geq_{Y_+}0$. Then, there exists $A \in \partial_{\epsilon'}^s\, G(\bar{x})$ 
but $A\not\in\partial_{\epsilon'+\epsilon}^e\,F(\bar{x})$. This will imply that\\
\begin{equation}\label{e1}
F(x^0) - F(\bar{x})- A(x^0-\bar{x})+\epsilon'+\epsilon\in -  Y_+\setminus l(Y_+)
\end{equation}\\[-.2cm]
for some $x^0\in X$ for which we also have\\
\begin{equation}\label{e2}
-G(x^0) + G(\bar{x}) + A(x^0-\bar{x}) - \epsilon' \in - Y_+.
\end{equation}\\[-.2cm]
Adding term by term (\ref{e1}) and (\ref{e2}) and taking into account that 
$-Y_+-Y_+\setminus l(Y_+)\subseteq -Y_+\setminus l(Y_+)$, we get \\[-.4cm]
$$
F(x^0) - G(x^0) \lneq_{Y_+} F(\bar{x}) - G(\bar{x}) - \epsilon,
$$\\[-.1cm]
i.e., $\bar{x}\not\in E_\epsilon^e (F-G,X,Y_+)$. This shows the direct implication.

\bigskip
\noindent
Sufficiency: we have by $({\cal H})$ that, for all $\epsilon' \gneq_{Y_+} 0$ 
and all $x\in \dom G$, there exists $A\in \partial_{\epsilon'}^s\, G(x)$. Hence, \\[-.2cm]
\begin{equation}\label{e3}
\forall x'\in X,\quad G(x') - G(x) \geq_{Y_+} A(x'- x) - \epsilon'.
\end{equation}\\[-.2cm]
Put $\epsilon''= G(\bar{x}) - G(x) - A(\bar{x}-x) + \epsilon'$ and substitute the value 
of $\epsilon'$ in (\ref{e3}) to immediately get \\[.1cm]
$$
\forall x'\in X,\quad G(x') - G(\bar{x}) \geq_{Y_+}  A(x'- \bar{x}) - \epsilon'',
$$\\[-.2cm]
with $\epsilon''\geq_{Y_+} 0$ because of (\ref{e3}) applied with $x'=\bar{x}$. Thus 
$A \in \partial_{\epsilon''}^s\, G(\bar{x})\subseteq\partial_{\epsilon''+\epsilon}^e\,F(\bar{x})$ 
by hypothesis. This implies that $\bar{x}\in\dom F\sqcap\dom G$ and that for all $x\in\dom G$,\\[-.1cm]
$$
F(x)-F(\bar{x})\not\not\lneq_{Y_+}A(x-\bar{x})-\epsilon''-\epsilon=G(x)-G(\bar{x})-\epsilon'-\epsilon, 
$$\\[-.3cm]
or equivalently,
\begin{equation}\label{e4}
F(x) - G(x) \not\not\lneq_{Y_+} F(\bar{x})  - G(\bar{x}) - (\epsilon' + \epsilon).
\end{equation}\\[-.2cm]
We can also assume that $x\in \dom F$ because otherwise, (\ref{e4}) obviously still 
holds. Since $\dom (F-G)=\dom F \sqcap \dom G$, it follows that (\ref{e4}) holds for 
any $x\in X$. This means together with $\bar{x}\in \dom F \sqcap \dom G=\dom (F-G)$ that 
$\bar{x}\in E_{\epsilon' + \epsilon}^e (F-G,X,Y_+)$, but $\epsilon' \gneq_{Y_+} 0$ was 
arbitrary, so according to Lemma \ref{reee}, $\bar{x}\in E_\epsilon^e (F-G,X,Y_+)$.\hfill$\Box$

\medskip 

\section{The main result}

Our main result in this paper establishes a first condition that completely characterizes the 
(nontrivial) $\varepsilon$-optimal solutions ($\varepsilon\geq 0$) to the following general (ROP): \\ 
$$
\varepsilon\text{\bf-}\Min_{h(x)\geq 0}\ f(x)
$$\\
where $f,\,h: X \longrightarrow \mathbb{R} \cup \{+\infty\}$ are such that $f$ is convex and 
$h$ is nonconcave. This extends and improves some results (e.g., \cite{hir,str,tse}) dealing 
with exact ($\varepsilon=0$) solutions.

\smallskip 

\begin{thm} \label{rop} Let $f,\,h: X \longrightarrow \mathbb{R} \cup \{+\infty\}$ 
be such that $f$ is convex and $h$ satisfies the following hypothesis: \\[-.3cm]
$$
({\cal H'})\qquad \forall\varepsilon \geq 0,\quad \forall x\in\dom h,\quad 
\partial_\varepsilon h(x) \neq \emptyset.
$$ \\[-.2cm]	
Let $\bar{x} \in \dom f$ such that $h(\bar{x})=0$, and, let $\varepsilon \geq 0$ such that  
$-\infty\leq \inf\limits_{x\in X}f(x) < f(\bar{x}) - \varepsilon$. Then, \\[-.2cm]
$$
\bar{x} \in \mathbf{\varepsilon}\text{\bf-}\argmin_{h(x)\geq 0} f(x)  \quad \Longleftrightarrow \quad 
\forall \varepsilon' \geq 0,\quad 
\partial_{\varepsilon'}h(\bar{x}) \subseteq \bigsqcup_{\alpha > 0} \partial_{\alpha\varepsilon+\varepsilon'}(\alpha f)(\bar{x}).
$$ 
\end{thm}

\smallskip 
\noindent 
{\em Proof.}\quad Necessity: suppose that  
$\bar{x} \in \mathbf{\varepsilon}\text{\bf-}\argmin_{h(x)\geq 0} f(x)$ and let 
$x^* \in \partial_{\varepsilon'}h(\bar{x})$ with $\varepsilon,\,\varepsilon' \geq 0$. 
Then, according to Lemma \ref{rb}, for $\epsilon=(\varepsilon,0)$, \\
$$
\bar{x} \in E^w_{\epsilon}\big((f,-h),X,\R^2_+\big) = 
E^w_{\epsilon}\big((f,0)-(0,h),X,\R^2_+\big).
$$ \\[-.2cm]
Following the hypothesis $({\cal H'})$ and using (\ref{ss}), the hypothesis $({\cal H})$ 
of Theorem \ref{cnswe} holds for $\sigma=w$ with $G=(0,h)$. In particular, for 
$\epsilon'=(0,\varepsilon')$, Theorem \ref{cnswe} applied to $F=(f,0)$ and $G=(0,h)$ 
shows that\\
$$
\partial^s_{\epsilon'}(0,h)(\bar{x})=\{0\} \times \partial_{\varepsilon'}h(\bar{x}) \; \subseteq \;  \partial^w_{\epsilon+\epsilon'}(f,0)(\bar{x}) = \partial^w_{(\varepsilon,\varepsilon')}(f,0)(\bar{x}).
$$ \\[-.1cm]
This implies that $(0,x^*) \in \partial^w_{(\varepsilon,\varepsilon')}(f,0)(\bar{x})$. 
Since the vector mapping $(f,0)$ is $\R^2_+$-convex, then by Theorem \ref{ssepw} of 
scalarization, there exists $(\lambda_1,\lambda_2) \in (\R^2_+)^w=\R^2_+\setminus\{0\}$ 
such that \\
\begin{equation}\label{e5}
\big< (\lambda_1,\lambda_2),(0,x^*)\big> = 
\lambda_2 x^* \in \partial_{\lambda_1 \varepsilon+\lambda_2 \varepsilon'}(\lambda_1 f)(\bar{x}).
\end{equation} \\[-.3cm]
Now, we will show that $\lambda_1$ and $\lambda_2$ are in fact both positive. Indeed, 
if $\lambda_1=0$ then $\lambda_2>0$ and (\ref{e5}) would imply that $x^*=0$. Since 
$x^* \in \partial_{\varepsilon'}h(\bar{x})$ was arbitrary, this would imply that 
$\partial_{\varepsilon'}h(\bar{x})=\{0\}$. But $\partial_{\varepsilon'}h(\bar{x})\supseteq 
\partial h(\bar{x})\neq \emptyset$ (by $({\cal H'})$) would imply that $0\in\partial h(\bar{x})$,  
meaning that $X=\{x : h(x)\geq 0\}$. Thus, with (\ref{eopta}), we would obtain that 
$\inf_{x\in X}f(x)=\inf_{h(x)\geq 0}f(x)\geq f(\bar{x})-\varepsilon$, which contradicts the 
assumption of the theorem. Likewise, if $\lambda_2=0$ then $\lambda_1 > 0$, 
and, (\ref{e5}) would imply that 
$0 \in \partial_{\lambda_1 \varepsilon}(\lambda_1 f)(\bar{x})=\lambda_1 \partial_{\varepsilon} f(\bar{x})$, 
which by definition means that $\bar{x} \in \mathbf{\varepsilon}\text{\bf-}\argmin_{x\in X} f(x)$, 
i.e., $\inf_{x\in X}f(x)\geq f(\bar{x})-\varepsilon$, which contradicts again the 
assumption of the theorem. Hence, $\lambda_1 > 0$ and $\lambda_2 > 0$. By setting 
$\alpha=\frac{\lambda_1}{\lambda_2}>0$ and substituting this value in (\ref{e5}), we finally get \\
$$
x^* \in \dfrac{1}{\lambda_2} \partial_{\lambda_1 \varepsilon+\lambda_2 \varepsilon'}(\lambda_1 f)(\bar{x}) = \partial_{\frac{\lambda_1}{\lambda_2} \varepsilon+\varepsilon'}\left(\frac{\lambda_1}{\lambda_2} f\right)(\bar{x}) = 
\partial_{\alpha \varepsilon+\varepsilon'}(\alpha f)(\bar{x}), 
$$ \\
which shows the necessary condition of the theorem. 

\bigskip
\noindent
Sufficiency: let 
$\epsilon=(\varepsilon,0)$, $\epsilon'=(\varepsilon'_1,\varepsilon'_2) \geq_{\R^2_+} 0$ 
and suppose that $z^*=(x^*,y^*) \in  \partial^s_{\epsilon'}(0,h)(\bar{x})$. Because 
$\partial^s_{\epsilon'}(0,h)(\bar{x}) = \{0\} \times \partial_{\varepsilon'_2}h(\bar{x})$, 
we have that $x^*=0$ and $y^* \in  \partial_{\varepsilon'_2}h(\bar{x})$. It follows, by the very 
hypothesis, that \\[-.1cm]   
$$
\exists \alpha > 0,\quad y^* \in 
\partial_{\alpha\varepsilon+\varepsilon'_2}(\alpha f)(\bar{x}) \subseteq 
\partial_{\alpha(\varepsilon+\varepsilon'_1)+\varepsilon'_2}(\alpha f)(\bar{x}), 
$$ \\[-.6cm]
or equivalently, \\[-.1cm]
$$
\exists \lambda=(\alpha,1) \in \Int \R_+^2, \quad \langle\lambda,z^*\rangle 
\in \partial_{\langle\lambda,\epsilon+\epsilon'\rangle}\big(\lambda \circ (f,0)\big)(\bar{x}).
$$  \\[-.2cm]
Since $(\R^2_+)^p=(\R^*_+)^\circ=\Int \R_+^2$, then by virtue of Theorem \ref{ssepw} 
and using (\ref{isdp}), it follows that 
$z^* \in \partial^{p}_{\epsilon+\epsilon'}(f,0)(\bar{x}) \subseteq 
\partial^{e}_{\epsilon+\epsilon'}(f,0)(\bar{x})$. But $z^* \in \partial^s_{\epsilon'}(0,h)(\bar{x})$ 
was arbitrary. This means that \\ 
$$ 
\forall \epsilon'\geq_{\mathbb{R}^{2}_{+}} 0,\quad \partial^s_{\epsilon'}(0,h)(\bar{x}) \subseteq \partial^{e}_{\epsilon+\epsilon'}(f,0)(\bar{x}).
$$ \\[-.2cm]
Thus, Theorem \ref{cnswe}, with $\sigma=e$, implies that \\ 
$$
\bar{x} \in E^e_{\epsilon}\big((f,0)-(0,h),X,\R^2_+\big) = E^e_{\epsilon}\big((f,-h),X,\R^2_+\big).
$$ \\[-.2cm] 
Using Lemma \ref{rb} and the assumption that $h(\bar{x})=0$, we finally get the desired 
result thus completing the proof of the theorem.\hfill$\Box$ 

\medskip 

\begin{rem} \label{rrp} \em Here are some remarks about the conditions of Theorem \ref{rop}: 
\begin{itemize}
\item[(a)] It has been proved in \cite[Corollary 2.1]{tuy} that if $f$ and $h$ are finite 
convex functions on $\R^n$, then the condition ``$h(\bar{x})=0$'' is necessary for the point 
$\bar{x}$ to be an exact solution for (ROP). Following our proof, this condition is only 
required to prove sufficiency, while the convexity of $f$ is required for necessity.
\item[(b)] The assumption ``$\inf_{x\in X}f(x) < f(\bar{x}) - \varepsilon$'' is only required 
to prove necessity. On the other hand, this assumption may be considered as ``essential'', since 
otherwise, we would have that $\bar{x}\in \mathbf{\varepsilon}\text{\bf-}\argmin_{x\in X} f(x)$, 
which obviously is characterized by $0\in \partial_{\varepsilon} f(\bar{x})$.
\item[(c)] The hypothesis $({\cal H'})$ is always fulfilled for all $\varepsilon > 0$ 
if, being given in a locally convex space, it holds that $h$ is proper convex lower 
semicontinuous on $X$ (see, e.g., \cite[Theorem 2.4.4.]{zali}). On the other hand, 
$({\cal H'})$ holds for $\varepsilon = 0$ if $h$ is proper convex lower semicontinuous 
on a Banach space $X$ and the Attouch--Brézis condition ``$\R_+[\dom h - x]$ is a closed 
vector subspace of $X$,'' holds for all $x\in\dom h$ (see \cite{att}). The Attouch--Brézis 
condition holds, for instance, if $x\in\Int\dom h$, in particular, if the convex function  
$h$ is continuous at $x$. So the convexity of $h$ is in some sense required in Theorem \ref{rop}.
\item[(d)] When $f,\,h : \R^n \longrightarrow \R$ are convex finite, hence continuous, 
so that $({\cal H'})$ is fulfilled, we recover the result by Hiriart-Urruty 
\cite[Theorem 3.5]{hir} characterizing the exact ($\varepsilon=0$) solutions to (ROP) 
(see also \cite[Remarks 3.6 and 3.7]{hir}): For $\bar{x}$ s.t. $h(\bar{x})=0$,\\
$$
\bar{x} \in \argmin_{h(x)\geq 0} f(x)  \quad \Longleftrightarrow \quad 
\forall \varepsilon \geq 0,\quad  
\partial_{\varepsilon}h(\bar{x}) \subseteq \bigsqcup_{\alpha > 0} \partial_{\varepsilon}(\alpha f)(\bar{x}).
$$\\
Note that in this case, the ``essential'' assumption reduces to 
$-\infty\leq \inf\limits_{x\in X}f(x) < f(\bar{x})$.
\end{itemize}
\end{rem}

\medskip

In the next section, two applications of Theorem \ref{rop} are derived for some 
important constrained cases. 	

\medskip

\section{Special reverse cases} 

\subsection{Additional convex constraints}

We first consider the special case where (ROP) is subject to 
additional convex constraints: \\[.1cm] 
$$
\varepsilon\text{\bf-}\Min_{h(x)\geq 0,\atop G(x)\leq 0}\ f(x)
$$\\[.1cm] 
where $f,\, h : X \longrightarrow \R\sqcup\{+\infty\}$ and  
$G=(g_1,\dots,g_m) : X \to \mathbb{R}^m\sqcup\{+\infty\}$ are convex mappings. 

\medskip  

In what follows, let $\Gamma(X)$ denote the set of proper convex functions 
$f : X \to \R\sqcup\{+\infty\}$, and let $\Gamma_0(X)$ denote the set of lower 
semicontinuous (l.s.c) functions $f$ in $\Gamma(X)$. \\

Using the penalty indicator function $\delta_C$ defined by: $\delta_C(x)=0$ 
if $x\in C$; $+\infty$ otherwise, we easily deduce the following $\varepsilon$-optimality 
criterion.  

\smallskip 

\begin{thm}\label{ropcc} 
Let $f,\, h : X \longrightarrow \R\sqcup\{+\infty\}$ and  
$G=(g_1,\dots,g_m) : X \to \mathbb{R}^m\sqcup\{+\infty\}$ be such that $h$ satisfies the 
hypothesis $({\cal H'})$, and, $(f,G)$ satisfy one of the following Moreau--Rockafellar 
$(M\!R)$ and Attouch--Brézis $(AB)$ qualification conditions:  \\
$$
\begin{array}{l}
(M\!R)\quad\left\{\begin{array}{l}
f\in\Gamma(X),\; g_i\in\Gamma(X)\;(\forall i),\; X \mbox{ separated locally convex},\\[.2cm]
\exists\, x_0\in \dom f \sqcap \dom G\ \mbox{ s.t. } 
\delta_{-\R_+^m} \mbox{ is continuous at }G(x_0).
\end{array}\right.\\[.6cm]
(AB)\quad\left\{\begin{array}{l}
f\in\Gamma_0(X),\; g_i\in\Gamma_0(X)\;(\forall i),\; X \mbox{ Fr\'echet space},\\[.2cm] 
\R_+[\R_+^m + G(\dom f\sqcap\dom G)] \mbox{ is a closed vector subspace of } \mathbb{R}^m.
\end{array}\right.
\end{array}
$$ \\[.1cm]
Let $\bar{x} \in \dom f$ s.t. $h(\bar{x})=0$, $G(\bar{x})\leq 0$, and, let $\varepsilon \geq 0$ 
s.t. $-\infty\leq \inf\limits_{G(x)\leq 0}f(x) < f(\bar{x}) - \varepsilon$. Then, \\[.1cm]
$$
\bar{x} \in \mathbf{\varepsilon}\text{\bf-}\argmin_{h(x)\geq 0,\atop G(x)\leq 0} f(x) 
\quad\Leftrightarrow\quad 
\forall \varepsilon' \geq 0, \quad \partial_{\varepsilon'}h(\bar{x}) \subseteq \!\!
\bigsqcup_{\alpha>0,\;\varepsilon_1\geq 0, \;\varepsilon_2\geq 0,\atop \varepsilon_1+\;\varepsilon_2=\alpha\varepsilon+\varepsilon'}\!\!
\bigsqcup_{\mu\in \R_+^m,\atop -\varepsilon_2 \leq \langle \mu,G(\bar{x})\rangle\leq 0}\!\!\!\!\!\!
\partial_{\varepsilon_1}(\alpha f + \mu\circ G)(\bar{x}).
$$ \\
In other words,
$$
\bar{x} \in \mathbf{\varepsilon}\text{\bf-}\argmin_{h(x)\geq 0,\atop G(x)\leq 0} f(x) 
\quad\Leftrightarrow\quad 
\left\{\begin{array}{l}
\forall \epsilon'\geq 0,\;\forall x^*\in \partial_{\epsilon'} h(\bar{x}),\;
\exists \alpha > 0,\; \exists \mu\in \R_+^m,\\[.3cm]
\exists\varepsilon_1\geq 0,\;\exists\varepsilon_2\geq 0,\; \varepsilon_1+\varepsilon_2=\alpha\varepsilon+\varepsilon' :\\[.3cm]
x^*\in \partial_{\varepsilon_1}(\alpha f + \mu\circ G)(\bar{x}),\\[.3cm]
-\varepsilon_2 \leq \langle \mu,G(\bar{x})\rangle\leq 0.
\end{array}
\right.
$$	
\end{thm}

\bigskip 
\noindent 
{\em Proof.}\quad The constrained (ROP) is penalized by the indicator function 
$\delta_{C}:=\delta_{-\R^m_+} \circ G$, where $C=\{x\in X : G(x)\leq 0\}$, so that \\
\begin{equation} \label{ropi}
\bar{x} \in \mathbf{\varepsilon}\text{\bf-}\argmin_{h(x)\geq 0,\atop G(x)\leq 0} f(x) 
\quad \Leftrightarrow \quad   
\bar{x} \in \mathbf{\varepsilon}\text{\bf-}\argmin_{h(x)\geq 0}\, \big(f + \delta_{-\R^m_+} \circ G\big)(x).
\end{equation} \\[-.1cm]
It is easy to check that Theorem \ref{rop} applies to the functions $f+\delta_C$ 
and $h$, that is to say that (\ref{ropi}) is equivalent to: 
\begin{equation}\label{ropic}
\forall \varepsilon' \geq 0, \quad \partial_{\varepsilon'}h(\bar{x}) \subseteq 
\bigsqcup_{\alpha > 0} 
\partial_{\alpha\varepsilon+\varepsilon'}\,\big(\alpha f+\alpha (\delta_{-\R^m_+} \circ G)\big)(\bar{x}).
\end{equation}   
According to the conditions $(M\!R)$ or $(AB)$ and the general hypotheses, by taking 
into account that $\alpha (\delta_{-\R^m_+} \circ G) = \delta_{-\R^m_+}\circ G$,    
the functions $f$, $\delta_{-\R^m_+}$ and $G$ satisfy exactly the qualification 
conditions of the Moreau--Rockafellar or Attouch--Brézis type and all the 
assumptions required for the $\varepsilon$-subdifferential addition/composition 
rule (see, e.g., \cite[Theorem 2.8.10, pp. 129]{zali}\footnote{The Attouch-Br\'ezis 
type condition (v) in \cite[Theorem 2.8.10]{zali} uses the li-convex concept for 
functions and sets which extends and englobes the convexity, l.s.c and closedness 
properties, see \cite[pp. 11, 15, 68]{zali}.}: \\[.1cm]
\begin{equation} \label{sdr}	
\partial_{\alpha\varepsilon+\varepsilon'}\,\big(\alpha f + \delta_{-\R^m_+} \circ G\big)(\bar{x}) = 
\bigsqcup_{\varepsilon_1\geq 0, \;\varepsilon_2\geq 0,\atop \varepsilon_1+\;\varepsilon_2=\alpha\varepsilon+\varepsilon'}
\;\bigsqcup_{\mu\in \partial_{\varepsilon_2}(\delta_{-\R^m_+})(G(\bar{x}))\;\sqcap\;\R^m_+}
\partial_{\varepsilon_1} (\alpha f + \mu\circ G)(\bar{x}).
\end{equation} \\[.1cm]
The conclusion of the theorem then follows from (\ref{ropi}), (\ref{ropic}), (\ref{sdr}) 
and the fact that, since $G(\bar{x})\leq 0$, \\[-.25cm]
$$
\partial_{\varepsilon_2}(\delta_{-\R^m_+})(G(\bar{x}))\; \sqcap\; \R^m_+ = 
\{\mu\in \R_+^m : \; -\varepsilon_2 \leq \langle \mu,G(\bar{x})\rangle\leq 0\}, 
$$ \\[-.2cm]
thus completing the proof of the theorem.\hfill$\Box$ 

\medskip

\begin{rem}\label{qc} \em
The Moreau--Rockafellar type condition $(M\!R)$ coincides in this application 
with the Slater condition ``$\exists\, x_0\in \dom f\sqcap\dom G\ \mbox{ s.t. }  
G(x_0)\in -\Int \R^m_+$,'' which is stronger in some sense than the 
Attouch--Brézis type condition $(AB)$. In fact, it easily implies that 
$\R_+[\R_+^m + G(\dom f\sqcap\dom G)]=\R^m$. But the Slater condition has 
the advantage of being in most cases easier to check. Also note that the 
reverse inclusion ``$\supseteq$'' in (\ref{sdr}) is easily proven without 
requiring any qualification conditions (see, e.g., \cite[Theorem 2.8.10]{zali} 
or \cite[Theorem 3.1]{elm2}).
\end{rem}

\subsection{Nonlinear equality constraint}

We derive in this section new $\varepsilon$-optimality criteria for the well-known 
nonlinear equality constrained optimization problem: 
$$
\varepsilon\text{\bf-}\Min_{h(x) = 0}\ f(x)
$$ 

This problem can first be viewed as (ROP) subject to an additional convex constraint 
as follows: 
$$
\varepsilon\text{\bf-}\Min_{h(x)\geq 0,\atop h(x)\leq 0}\ f(x)
$$ 

To the best of our knowledge, the following result has not been previously established 
by specialists even for dealing with exact solutions. 

\smallskip 

\begin{cor}\label{ropce} 
Let $f,\, h : X \longrightarrow \R\sqcup\{+\infty\}$ be convex proper on $X$ a separated 
locally convex space such that $h$ satisfies the hypothesis $({\cal H'})$. 
Let $\bar{x} \in \dom f$ s.t. $h(\bar{x}) = 0$, and, let $\varepsilon \geq 0$ s.t.  
$-\infty \leq \inf\limits_{h(x)\leq 0}f(x) < f(\bar{x}) - \varepsilon$. Then, \\
$$
\bar{x} \in \mathbf{\varepsilon}\text{\bf-}\argmin_{h(x) = 0} f(x) 
\quad\Leftrightarrow\quad 
\forall \varepsilon' \geq 0, \quad \partial_{\varepsilon'}h(\bar{x}) \subseteq 
\displaystyle\bigsqcup_{\alpha > 0,\;\beta\geq 0} \partial_{\alpha\varepsilon+\varepsilon'}(\alpha f+\beta h)(\bar{x}).
$$
\end{cor}

\medskip 
\noindent 
{\em Proof.}\quad We first check that all the assumptions of Theorem \ref{ropcc} 
are satisfied for $f$, $h$ and $G=h$. Indeed, following Remark \ref{qc}, the 
Moreau--Rockafellar type condition $(M\!R)$ which coincides in this case with the 
Slater condition ``$\exists\, x_0\in \dom f\sqcap\dom h\ \mbox{ s.t. } h(x_0) < 0$,'' 
is well guaranteed under the ``essential'' hypothesis of the corollary. Indeed, if 
$h(x)\geq 0$ for all $x\in\dom f\sqcap\dom h$, we would have  
$\{x \in \dom f : h(x) \leq 0\}=\{x \in \dom f : h(x) = 0\}$. Thus, for 
$\bar{x} \in \mathbf{\varepsilon}\text{\bf-}\argmin_{h(x) = 0} f(x)$, we would obtain 
that $\inf_{h(x)\leq 0}f(x)=\inf_{h(x) = 0}f(x) \geq f(\bar{x})-\varepsilon$, which 
contradicts the ``essential'' hypothesis that $\inf_{h(x)\leq 0}f(x) < f(\bar{x})-\varepsilon$. 
Hence, for the direct implication, applying Theorem \ref{ropcc}, we get for all 
$\varepsilon'\geq 0$, \\ 
$$
\begin{array}{lcl}
\partial_{\varepsilon'}h(\bar{x}) 
& \subseteq &
\displaystyle \bigsqcup_{\alpha>0,\;\varepsilon_1\geq 0, \;\varepsilon_2\geq 0,\atop \varepsilon_1+\;\varepsilon_2=\alpha\varepsilon+\varepsilon'}
\bigsqcup_{\mu\in \R_+,\atop -\varepsilon_2 \leq \langle \mu,h(\bar{x})\rangle\leq 0}
\partial_{\varepsilon_1} (\alpha f + \mu h)(\bar{x}) \\[1cm]
& = & \displaystyle \bigsqcup_{\alpha>0,\;\;\beta\geq 0,\atop 0\leq\varepsilon_1\leq\alpha\varepsilon+\varepsilon'} 
\partial_{\varepsilon_1} (\alpha f + \beta h)(\bar{x}) \\[1cm]
& = & \displaystyle \bigsqcup_{\alpha>0,\;\beta\geq 0} 
\partial_{\alpha\varepsilon+\varepsilon'} (\alpha f + \beta h)(\bar{x}).
\end{array} 
$$ \\[-.1cm]
The reverse implication follows easily by the same arguments as for Theorem \ref{ropcc} 
by taking into account Remark \ref{qc} about the inclusion ``$\supseteq$'' for  
the formula (\ref{sdr}).\hfill$\Box$  

\bigskip

Let us recall the following results (see  also Remark \ref{rrp}(a)) that hold in the 
particular case of convex finite functions. This justifies why the solutions of 
(ROP) are taken on the boundary of the reverse constraint set. 

\smallskip 

\begin{lem} {\em (\cite{tuy})} \label{cccf} 
Let $f,\, h : X \to \R$ be finite convex functions and the essential assumption 
``$-\infty \leq \inf\limits_{x\in X}f(x) < \inf\limits_{h(x)\geq 0}f(x)$'' is 
satisfied. Then, for every $x\in X$ such that $h(x)>0$, there exists $\pi(x)$ 
such that $h(\pi(x))=0$ and $f(\pi(x)) < f(x)$. In other words, \\[-.1cm]
$$
\argmin_{h(x) = 0} f(x) = \argmin_{h(x) \geq 0} f(x).
$$
\end{lem}

\bigskip

We similarly derive a result for approximate solutions that holds true 
under the same conditions.

\smallskip 

\begin{cor} \label{ecccf}
Under the conditions of Lemma \ref{cccf},  the following property holds: 
$\forall \varepsilon\geq 0$, \\[-.1cm]
$$
\mathbf{\varepsilon}\text{\bf-}\argmin_{h(x)= 0} f(x) = 
\mathbf{\varepsilon}\text{\bf-}\argmin_{h(x)\geq 0} f(x) \sqcap \{x\in X :\; h(x) = 0\}.
$$
\end{cor}

\bigskip 
\noindent 
{\em Proof.}\quad Let us prove the first inclusion. By the definition (\ref{eopta}), 
$\bar{x} \in \mathbf{\varepsilon}\text{\bf-}\argmin_{h(x)= 0} f(x)$, iff
$$
h(\bar{x})= 0, \qquad  f(\bar{x}) \leq \inf_{h(x)= 0} f(x) + \varepsilon.
$$
Let $x'\in X$ such that $h(x')\geq 0$. If $h(x')=0$ then 
$f(\bar{x}) \leq \inf_{h(x)= 0} f(x) + \varepsilon \leq f(x') + \varepsilon$. 
If $h(x') > 0$ then, by Lemma \ref{cccf}, there exists $\pi(x')$ 
s.t. $h(\pi(x'))=0$ and $f(\pi(x')) < f(x')$. Hence, we also obtain that  
$f(\bar{x}) \leq f(\pi(x')) + \varepsilon < f(x') + \varepsilon$. So, in all cases,   
$$
h(\bar{x}) = 0, \qquad f(\bar{x}) \leq \inf_{h(x) \geq 0} f(x) + \varepsilon,
$$ 
i.e., $\bar{x} \in \mathbf{\varepsilon}\text{\bf-}\argmin_{h(x)\geq 0} f(x)$. 
Conversely, let  $\bar{x} \in \mathbf{\varepsilon}\text{\bf-}\argmin_{h(x) \geq 0} f(x)$ 
such that $h(\bar{x}) = 0$. Then, it is immediate that 
$$
f(\bar{x}) \leq \inf_{h(x) \geq 0} f(x) + \varepsilon \leq \inf_{h(x)= 0} f(x) + \varepsilon, 
$$
which shows that $\bar{x} \in \mathbf{\varepsilon}\text{\bf-}\argmin_{h(x)= 0} f(x)$.\hfill$\Box$

\bigskip 

Under the statements of Corollary \ref{ecccf}, we straightforwardly derive 
from Theorem \ref{rop} the following simplified $\varepsilon$-optimality 
criterion for nonlinear equality constrained convex optimization problems. 
  
\smallskip 
    
\begin{cor} \label{sropec}
Under the assumptions of Lemma \ref{cccf}, let $\bar{x} \in \dom f$ such that $h(\bar{x}) = 0$, 
and, let $\varepsilon \geq 0$ such that $-\infty \leq \inf\limits_{h(x)\leq 0}f(x) < f(\bar{x}) - \varepsilon$. 
Then,  \\[.1cm]
$$
\bar{x} \in \mathbf{\varepsilon}\text{\bf-}\argmin_{h(x)=0} f(x)  \quad \Longleftrightarrow \quad 
\forall \varepsilon' \geq 0,\quad \partial_{\varepsilon'}h(\bar{x}) \subseteq 
\bigsqcup_{\alpha > 0} \partial_{\alpha\varepsilon+\varepsilon'}(\alpha f)(\bar{x}).
$$ 
\end{cor}

\smallskip 

\begin{rem} \em Here are some remarks about the conditions of Corollary \ref{sropec}: 
\begin{itemize}
\item[(a)] The assumption ``$\inf_{h(x)\leq 0}f(x) < f(\bar{x}) - \varepsilon$'' can also be 
considered here as essential, since otherwise, we would have that 
$\bar{x}\in \mathbf{\varepsilon}\text{\bf-}\argmin_{h(x)\leq 0} f(x)$, which is a trivial case of 
convex optimization that may be easily characterized by the following condition: 
$$
0\in \bigsqcup_{\beta\geq 0} \partial_{\varepsilon} (f + \beta h)(\bar{x})
$$ 
\item[(b)] Following Remark \ref{rrp}(c), the hypothesis $({\cal H'})$ of Theorem \ref{rop} 
is always fulfilled under the assumptions of Lemma \ref{cccf}. However, the 
restrictive hypothesis on the functions $f$ and $h$, which must be finite 
throughout the space $X$, does not allow to directly extend the previous 
result to nonlinear equality constrained convex optimization problems subject 
to additional convex constraints. Thus, this problem remains open for establishing 
a similar result that holds true even for extended valued convex mappings, so that 
Theorem \ref{ropcc} may be applied instead of Theorem \ref{rop}.
\end{itemize}
\end{rem} 

\medskip

\noindent
{\footnotesize {\bf Statements and Declarations.} 
The authors have no pertinent declarations concerning conflicts of interest, financial or 
non-financial interests, competing interests or other statements to disclose.}

\end{document}